# Counterexamples in the Cutting Stock Problem

Constantine Goulimis

Greycon Ltd., 7 Calico House, Plantation Wharf, London SW11 3TN.
cng@greycon.com

**Abstract**

We present a collection of examples in various aspects of the one-dimensional cutting stock problem ('1D-CSP'). Based on our 37+ year experience in this area, we discuss some modelling pitfalls and counterexamples.

Keywords: Cutting stock problem, pattern reduction, setup minimisation, integer roundup, paper industry, plastic film industry

## 1 Introduction

The many practical applications of the one-dimensional cutting stock problem ('1D-CSP') have provided a rich source of challenges to the mathematical optimisation community.

We use the following terminology for the 1D-CSP:

- The master size is $W$

- The orders, indexed over $j$ ($1 \leq j \leq m$), have size $w_j$ and we must produce in the range $[q_j, Q_j]$. We assume, without loss of generality, that the sizes are ordered $w_1 > w_2 > \ldots > 0$.

- We associate with each possible pattern $i$ (defined by the vector $a_i = (a_{ij})$ of the number of times order $j$ appears in it) an integer variable $x_i$ (= the number of repetitions of this pattern in the solution) and a waste value of $c_i$.

The *general* formulation is:

**min** $\quad z = \sum_i c_i x_i$

**subject to** $\quad Q_j \geq \sum_i a_{ij} x_i \geq q_j, \forall j$

In terms of the objective function, we distinguish between:

- *Minimisation of the number of master items*, when $c_i = 1$.

- *Waste minimisation*, when $c_i$ is the waste associated with pattern $i$.

In terms of the quantity constraints on the orders, there are two special cases:

- *One-sided*, when we only have a lower bound for the order quantity: $\sum_i a_{ij} x_i \geq q_j, \forall j$

- *Equality-constrained*, when $Q_j = q_j, \forall j$.

In this paper we present counterexamples that shine a light on three areas:

1. Master item vs. waste minimisation, which is more appropriate?
2. Order splitting, i.e. is there a limit to how many patterns contain an order?
3. Pattern count, what can we tell about the number of patterns in an optimal solution?

We discuss these in the following sections. Our work is based on our experience, which started with the first-published integer optimal algorithm to the 1D-CSP [1] and continued over the last 37 years with the commercial delivery of a software tool (X-Trim) in use mostly in the paper, plastic film and metals industries in 41 countries. In this paper, we discuss some of the proven and unproven conjectures and assumptions. We also propose a new conjecture.

## 2 Integer Roundup

If the objective is to minimise the number of master items, one conjecture relates to the value of the linear programming relaxation. Marcotte in 1985 [2] considered the traditional one-sided / minimising the number of master items formulation:

**min** $\quad z = \sum_i x_i$

**subject to** $\quad \sum_i a_{ij} x_i \geq q_j, \forall j$

and claimed that the integer round-up property:

$$z^* = \lceil z_{LP} \rceil$$

would hold (i.e. the linear programming result, when rounded up, would be the same as the integer optimal solution). So, if the linear programming relaxation produces a solution with 12.3 master items, then, the theory goes, the optimal (minimum number of master items) solution would require $\lceil 12.3 \rceil$ = 13. She proved this property for certain specific classes of the 1D-CSP, notably those where there are only two orders and those where the sizes have successive divisibility, $w_1 \mid w_2 \mid \ldots$ . These are very special cases indeed.

In a subsequent paper [3], Marcotte found some counterexamples to the above rule (with a gap of exactly 1), conjecturing that these counterexamples would require coefficients of magnitude $10^7$. She also wrote "This observation supports the claim made in Section 1, that any counterexample to the rounding property is unlikely to arise from a practical problem".

This turned out [4] not be true either, the following instance

$W = 132 \quad\quad \mathbf{w} = (44,33,12) \quad \mathbf{q} = (2,3,6)$

has $z_{LP}$ = 259/132 = 1.962… whereas $z^*$ = 3, therefore the gap is

$$\frac{137}{132} = 1.0378 \ldots$$

Yet the myth of integer roundup would persist; papers even in 1993 [5] would state "integer round-up property seems to hold for almost any real-world problem".

In 2020 Ripatti & Kartak [6] managed to prove that the property does hold, for all equality-constrained instances where $W \leq 15$. This is unfortunately too small to be of practical interest.

Lastly, Scheithauer & Terno [7], moving in a different direction, conjectured that there is a *modified integer round-up property* that applies to this formulation of the CSP, namely a gap of at most 2:

$$z^* \leq 1 + \lceil z_{LP} \rceil$$

This conjecture has not, to our knowledge, been proven, although we also suspect it is true. Assuming it to be true, Scheithauer & Terno went on to propose algorithms that exploit it.

However, even if the conjecture were true, the underlying formulation is irrelevant for most practical purposes, because minimising the number of master items is not the same as minimising waste (as we shall discuss in the next section).

## 3 Minimising Master Items is Not the Same as Minimising Waste

Many of the research papers that have been published, only consider lower limits on the order quantities, as opposed to the two-sided constraints of the general formulation. Equally, many papers consider the special case of minimising the number of master items instead of waste minimisation.

In this section we demonstrate why neither is a good idea. We do this by showing that

- Minimising the number of master items does not minimise waste, and
- Minimising waste does not minimise the number of master items.

Since this implies that the two objectives have different outcomes (except in the equality-constrained special case), the question becomes which one of the two is more important in practice? In our experience, waste minimisation is a better objective in practically all of the situations we have encountered.

### 3.1 Minimum Number of Master Items Does Not Imply Minimum Waste

We start with obvious statement that for equality-constrained problems (no tolerance on the order quantity) waste minimisation and master item minimisation are equivalent.

But for the general problem (or even the one-sided special case), does a solution that minimises the number of master items, also minimise waste? We start with the most trivial example that demonstrates that this is not so. Suppose we have single order that fits $a$ ($>2$) times in the master item and that $q = a+1$, i.e. we require at a minimum $a+1$ copies of this order.

Then the two-pattern solutions containing $a+1$, $a+2$, …, $2a-1$ in total all require 2 master items, but the waste differs massively. Here are a couple of solutions when $a = 3$, both of whom require the same number of master items:

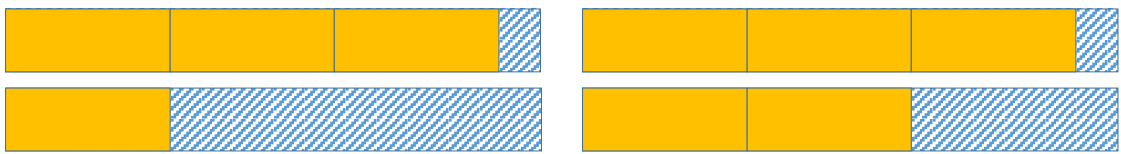

This counterexample however breaks down in one-sided scenarios, if we demand, as is common, that patterns are *maximal*, i.e. the wasted space is $< \min(w_i)$. So we need to move to more complicated counterexamples.

Our first counterexample comes from a paper by Cui et al [8], who present the following solution, generated by their proposed algorithm, to an instance (redrawn to make it clearer):

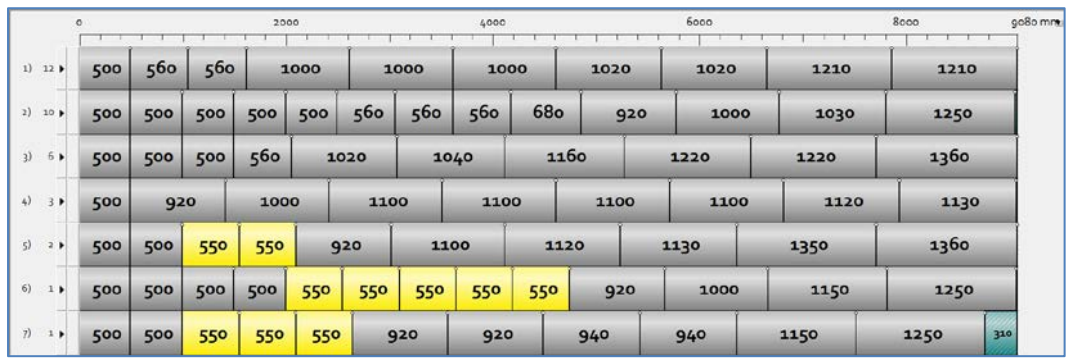

The above solution has waste ≈ 0.173%, consumes 35 master items and produces 12 copies of size 550 (highlighted). Cui et al use the one-sided formulation.

In fact, solutions with the same number (35) of master items, 7 patterns and 0.0% waste are possible:

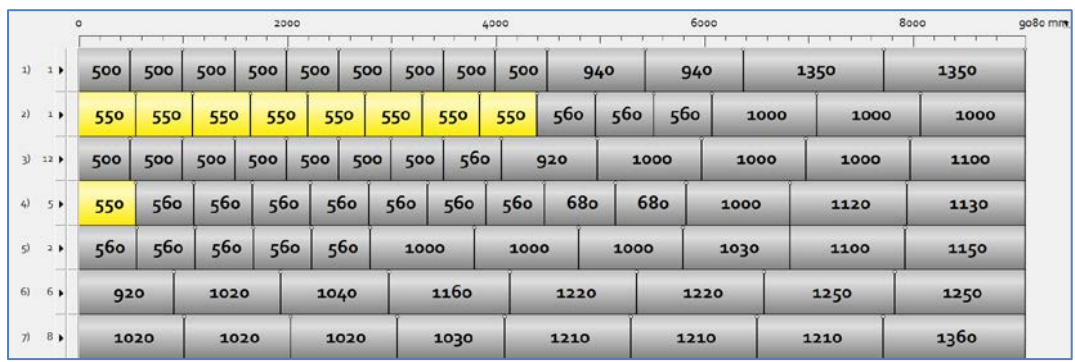

The above solution, which is feasible in the one-sided formulation, contains one more item (12→13) of size 550, within the same master item and pattern count as the solution in the paper.

So, in practical terms, minimising the number of master items in the presence of one-sided constraints can lead to solutions with excess waste.

This generalises trivially to the general two-sided case, all we have to do is constrain, e.g., the item of size 550 to appear [12, 13] times.

In addition to what may have been theoretical papers, we have often encountered this scenario in practice with X-Trim. Here is a real-world example with 4 orders and 63 master items against a competitive product:

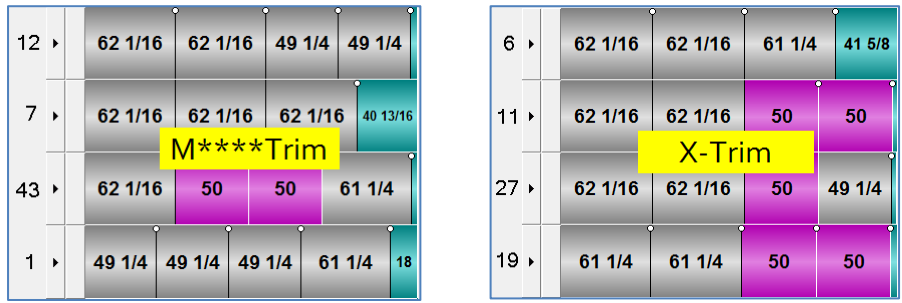

The X-Trim solution contains an additional item of size 50 (87 *vs*. 86).

### 3.2 Waste Minimum Does Not Imply Master Item Minimum

The question then arises: does the implication hold the other way round? I.e. does minimising waste also minimise the number of master items? Again, we have to exclude the trivial equality-constrained case. In this section we provide a counterexample that not only demonstrates that this is not true, but also that the difference in the number of master items between the two types of solutions is not small.

We start with two-sided constraints:

$$W = 1000 \qquad w = (300, 340) \qquad q = (15, 15) \qquad Q = (15, 18)$$

Then the minimum-waste solution contains 11 master items. However, the minimum-master-items solution consumes 10, but has higher waste:

|  | Waste-minimal Solution | Master-item-minimal Solution |
|---|---|---|
| Solution | | |
| Waste | 380 = 7×20 + 4×60 | 400 = 5×20 + 5×60 |
| Waste | 3.455% | 4.000% |
| Master | 11 | 10 |

Having shown by this counterexample that a minimum waste solution may actually require more than the minimum number of master items for two-sided problems, we might imagine that for the special (even if not practically interesting) case of one-sided constraints the conjecture might be true. Alas, the above example can be extended to this case also. Varying the above counterexample by removing the upper bound $Q$ on the orders:

$$W = 1000 \qquad w = (300, 340) \qquad q = (15, 15) \qquad Q = (\infty, \infty)$$

we obtain the following waste-optimal solution:

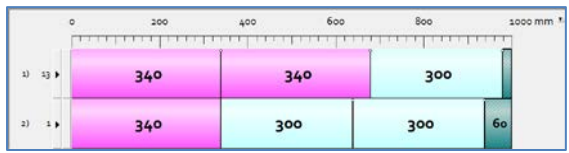

This however uses 14 master items (whereas a minimum-master-item solution, as we have seen above, uses only 10). Not only the minimum-master-item solution has fewer master items than the minimum-waste solution, but the gap is quite large: in this instance, the minimum waste solution (with 14 master items) has 40% more master items than the minimum-master-items one (with 10).

In conclusion, waste minimisation does not imply minimisation of master items, except in the equality constrained case.

# 4 Order Splitting

In a particular solution, an order is *split*, if it appears in more than one pattern. Splitting is undesirable because it may lead to difficulties in shipping the production.

Is there some upper bound on the number of times an order is split in an optimal solution? It is tempting to assume that splitting can be controlled. So much so, that at least one algorithm has been proposed using the assumption that every order will appear in at most 2 patterns.

In this section we provide a counterexample to this premise.

## 4.1 Order Appears in $k$ Patterns

Consider, again, the trivial case where a single order fits $a$ ($= \lfloor W/w_1 \rfloor > 2$) times in the master item and we require (exactly) $Q$ copies. If $Q \neq 0 \bmod a$, the unique optimal (both in terms of waste and in terms of number of master items) solution has two patterns. The first pattern contains $a$ times the single order and will be reproduced $\lfloor Q/a \rfloor$ times. The second pattern will consume one master item and will contain the balance (if any). So, for this rather trivial scenario, we have one order appearing in two patterns. So, for equality-constrained problems at least, an order may have to appear in at least two patterns.

But it gets much worse: in the more general case, in the presence of multiple orders, an order may **have** to appear in $k$ patterns for arbitrary $k > 0$. Here is a trivial example with $k = 3$:

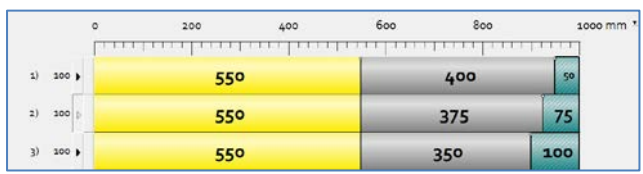

This example is even relevant to the special cases (minimise number of master items, one-sided constraints and equality constraints).

This counterexample is interesting because it refutes an approach by Haessler & Talbot [9] for solving corrugator trim optimisation problems using set partitioning. In this approach, a solution, comprises of *elements* of four types:

    a. An order is (fully) satisfied by one pattern
    b. Two orders (fully) satisfied by a single pattern
    c. Two orders (fully) satisfied by two patterns
    d. Three orders (fully) satisfied by two patterns

In Haessler & Talbot's approach each element has a waste level and then the optimisation problem becomes a question of solving a set partitioning over these elements.

However, in this approach an order can appear in at most two patterns. Hence this approach is sub-optimal (in fact it would falsely declare the above counterexample to be infeasible).

This counterexample depends crucially on the fact that the largest order has $w_1$ is in the range ($W/2$, $2W/3$) and all the others have size in the range ($W/3$, $W/2$). What happens if all orders have size < $W/2$ is not currently known, but the following $m = 3$ variant shows that we can have all widths smaller than $W/2$, yet one order has to appear in every single pattern:

$$W = 1000 \qquad w = (400, 375, 350, 275) \qquad q = Q = (1, 1, 1, 300)$$

whose optimal solution requires that the 275 size appears in every pattern:

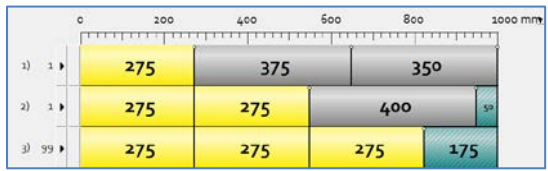

# 5 Pattern Count

The 1D-CSP is quite degenerate in the sense that there are often many different waste-optimal solutions. Within this universe of waste-optimal solutions, it is quite important in some industrial applications to find the solution(s) with the minimum pattern count. This is a notoriously hard problem [10] [8] [11] [12].

Linear programming theory suggests that a 1D-CSP with $m$ orders will require the same number of patterns. This is no more than a suggestion, because

- We are dealing with integer programming, not linear
- There may be degeneracy in the linear programming basis

In this section we provide examples on the number of patterns in the range [1, ($m+2$)] for equality-constrained instances with $m$ orders. It is interesting here to compare with the average behaviour: empirical evidence is that we need ~0.6 $m$ – 0.7 $m$ patterns on average for $m$ orders [13].

## 5.1 Equality-constrained - $m$ Orders Require ($m+2$) Patterns

Some authors, e.g. [14], have mistakenly assumed that with $m$ orders, $m$ patterns are sufficient. Instances where $m$ orders require $m+1$ patterns for equality-constrained problems are easy to construct. Starting with a trivial one-order, equality-constrained, problem that requires two patterns:

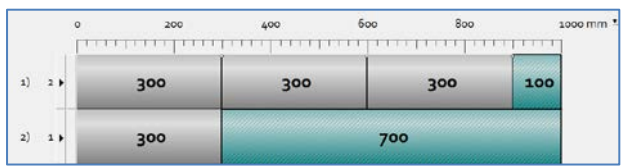

we can extend this to multiple orders, e.g. here is an example with $m = 2$:

$$W = 1000 \qquad w = (300, 250) \qquad q = Q = (7, 100)$$

The 3-pattern optimal solution is this:

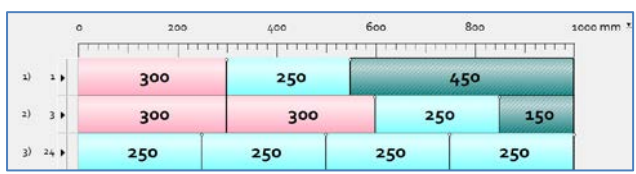

To understand how rare instances with $m$ orders and ($m+1$) patterns are, we generated 1,000 (realistic) instances with 10 orders each. Of these, 2 instances required 10 patterns, 2 instances required 11 and all the rest (i.e. 99.6%) required 9 or less.

For many years the author had therefore believed that the following conjecture is true:

> Conjecture A: for an equality-constrained problem with $m$ orders, $(m+1)$ patterns are sufficient, i.e. there is a minimum waste solution with no more than $(m+1)$ patterns.

However, whilst conducting similar computational studies, we came across a 9-order instance whose optimal solution requires 11 patterns:

$W = 560$

$w = (109, 114, 115, 125, 132, 140, 142, 148, 200)$

$Q = q = (46, 39, 39, 36, 22, 23, 32, 22, 42)$

We have established computationally that there are no (waste optimal) solutions to this instance with fewer than 11 patterns; one such solution is shown below:

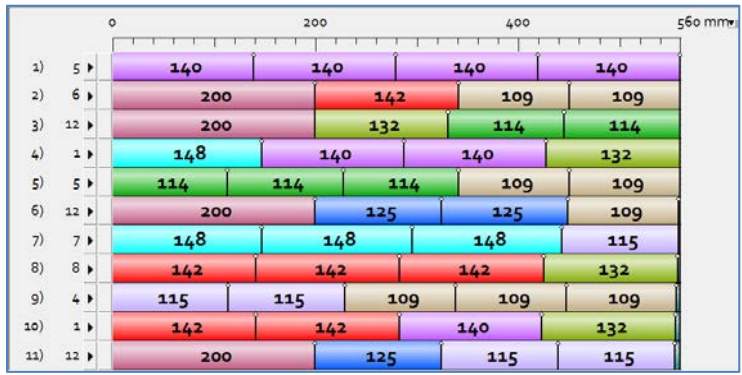

Thus Conjecture A does not hold! This is the simplest case we have discovered to date where $(m+2)$ patterns are needed.

Perhaps therefore the following conjecture is true:

> Conjecture B: for an equality-constrained problem with $m$ orders, $(m+2)$ patterns are sufficient, i.e. there is a minimum waste solution with no more than $(m+2)$ patterns.

### 5.2 Equality-constrained - $m$ Orders Require One Pattern

At the other end of the spectrum, artificial instances with (arbitrary) $m$ orders requiring only one pattern are also possible. Here is an example with $m = 5$,

$W = 1000$    $w = (300, 250, 200, 150, 90)$    $q = Q = (100, 100, 100, 100, 100)$

The unique optimal solution is:

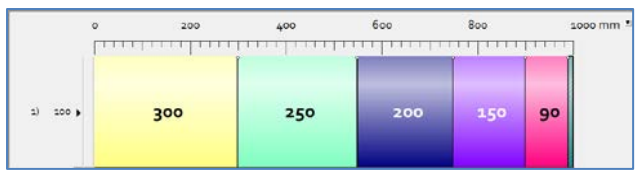

Note that this applies to all special cases of the 1D-CSP.

Interestingly, this construction gives rise to the only easily-calculated lower bound on the number of patterns. The minimum number of patterns required is:

$$\left\lceil \frac{\sum_j w_j}{W} \right\rceil$$

In practice, this sharp lower bound (which does not depend on $q$, $Q$) is very poor, except for very small instances.

We can generalise the above construction so that $m$ orders require $k \leq m$ patterns in a straightforward manner.

### 5.3 Implications for the General Formulation

If Conjecture B is true, it implies that any instance of 1D-CSP, regardless of the objective or the nature of the bounds, has an optimal solution with no more than ($m$+2) patterns. The proof is by contradiction: if there were an instance arising from the general formulation that did require more than ($m$+2) patterns, then the corresponding equality-constrained instance would violate the conjecture.

The converse is not true: it is possible e.g. that a stronger bound exists for one-sided problems, perhaps $m$, or, possibly even below this.

Practical instances very often have side constraints that make certain patterns infeasible (e.g. constraints on the number of slitting knives). In the presence of such constraints, it is probably impossible to make statements on the pattern count.

# 6 Acknowledgments

The author acknowledges the helpful comments by Sophia Drossopoulou.